\input amstex\documentstyle {amsppt}  
\pagewidth{12.5 cm}\pageheight{19 cm}\magnification\magstep1
\topmatter
\title Remarks on Springer's representations\endtitle
\author G. Lusztig\endauthor
\address Department of Mathematics, M.I.T., Cambridge, MA 02139\endaddress
\thanks Supported in part by the National Science Foundation\endthanks
\endtopmatter   
\document

\define\hcx{\hat{\cx}}

\define\Irr{\text{\rm Irr}}

\define\sgn{\text{\rm sgn}}

\define\ua{\un a}

\define\uc{\un c}

\define\da{\dagger}

\define\Lie{\text{\rm Lie }}

\define\frl{\forall}

\define\sqc{\sqcup}

\define\qua{\quad}

\define\part{\partial}

\define\n{\notin}

\define\m{\mapsto}
\define\do{\dots}

\define\lra{\leftrightarrow}

\define\sub{\subset}    
\define\bxt{\boxtimes}
\define\T{\times}
\define\ti{\tilde}
\define\nl{\newline}
\redefine\i{^{-1}}

\define\un{\underline}

\define\ot{\otimes}
\define\bbq{\bar{\QQ}_l}

\define\Hom{\text{\rm Hom}}
\define\End{\text{\rm End}}

\define\Ind{\text{\rm Ind}}

\define\card{\text{\rm card}}

\define\a{\alpha}
\redefine\b{\beta}

\define\p{\pi}
\define\ph{\phi}

\define\r{\rho}

\redefine\t{\tau}

\define\z{\zeta}
\define\x{\xi}

\define\Th{\Theta}

\define\kk{\bold k}

\define\NN{\bold N}

\define\QQ{\bold Q}

\define\TT{\bold T}

\define\WW{\bold W}
\define\ZZ{\bold Z}

\define\YY{\bold Y}

\define\ca{\Cal A}
\define\cb{\Cal B}
\define\cc{\Cal C}
\define\cd{\Cal D}
\define\ce{\Cal E}

\define\ch{\Cal H}

\define\cl{\Cal L}

\define\cn{\Cal N}
\define\co{\Cal O}

\define\cs{\Cal S}
\define\ct{\Cal T}
\define\cu{\Cal U}

\define\cx{\Cal X}

\define\fg{\frak g}

\define\ft{\frak t}

\define\fS{\frak S}

\define\ta{\ti a}
\define\tb{\ti b}
\define\tc{\ti c}

\define\tE{\ti E}

\define\tcb{\ti{\cb}}

\define\ALV{A}
\define\HS{HS}
\define\LCL{L1}
\define\LUG{L2}
\define\ORA{L3}
\define\LUI{L4}
\define\LUU{L5}
\define\LUS{L6}
\define\LS{LS}
\define\SE{Se}
\define\SPA{S1}
\define\SPAS{S2}
\define\SPAF{S3}
\define\SPR{Sp}
\define\XUE{X}

\head Introduction\endhead
\subhead 0.1\endsubhead
Let $\kk$ be an algebraically closed field of characteristic exponent $p\ge1$. Let $G$ be a connected reductive 
algebraic group over $\kk$ and let $\fg$ be the Lie algebra of $G$. Let $\cu_G$ be the variety of unipotent 
elements of $G$ and let $\cn_\fg$ be the variety of nilpotent elements of $\fg$ (we say that $x\in\fg$ is 
nilpotent if for some/any closed imbedding $G\sub GL(\kk^n)$, the image of $x$ under the induced map of Lie
algebras $\fg@>>>\End(\kk^n)$ is nilpotent as an endomorphism). Note that $G$ acts on $G$ and $\fg$ by the adjoint
action. Let $\cx_G$ (resp. $\cx_\fg$) be the set of $G$-orbits on $\cu_G$ (resp. on 
$\cn_\fg$). We fix a prime number $l$, $l\ne p$. Let $\hcx_G$ (resp. $\hcx_\fg$) be the set of pairs $(\co,\cl)$ 
where $\co\in\cx_G$ (resp. $\co\in\cx_\fg$) and $\cl$ is an irreducible $G$-equivariant $\bbq$-local system on 
$\co$ up to isomorphism. Let $\WW$ be the Weyl group of $G$. For any Weyl group $W$ let $\Irr(W)$ be the set of 
isomorphism classes of irreducible representations of $W$ over $\QQ$. In \cite{\SPR}, Springer defined (assuming 
that $p=1$ or $p\gg0$) natural injective maps $S_G:\Irr(\WW)@>>>\hcx_G$, $S_\fg:\Irr(\WW)@>>>\hcx_\fg$ (each of 
these two maps determines the other since in this case we have canonically $\hcx_G=\hcx_\fg$). In \cite{\LUG} a
new definition of the map $S_G$ (based on intersection homology) was given which applies without restriction on 
$p$. A similar method can be used to define $S_\fg$ without restriction on $p$ (see \cite{\XUE} and 2.2 below); 
note that in general $\hcx_G,\hcx_\fg$ cannot be identified. Now for any $\co\in\cx_G$ (resp. $\co\in\cx_\fg$),
$(\co,\bbq)$ is in the image of $S_G$ (resp. $S_\fg$) hence there is a well defined injective map 
$S'_G:\cx_G@>>>\Irr(\WW)$ (resp. $S'_\fg:\cx_\fg@>>>\Irr(\WW)$) such that for any $\co\in\cx_G$ (resp. 
$\co\in\cx_\fg$) we have $S'_G(\co)=E$ (resp. $S'_\fg(\co)=E$) where $E\in\Irr(\WW)$ is given by
$S_G(E)=(\co,\bbq)$ (resp. $S_\fg(E)=(\co,\bbq)$). Let $\fS_G$ be the image of $S'_G:\cx_G@>>>\Irr(\WW)$. Let 
$\fS_\fg$ be the image of $S'_\fg:\cx_\fg@>>>\Irr(\WW)$. 

In \cite{\LUU}, we gave an apriori definition (in the framework of Weyl groups) of the subset $\fS_G$ of 
$\Irr(\WW)$ which parametrizes the unipotent $G$-orbits in $G$. In this paper we give an apriori definition (in a
similar spirit) of the subset $\fS_\fg$ of $\Irr(\WW)$ which parametrizes the nilpotent $G$-orbits in $\fg$. (See
Proposition 3.2.) This relies heavily on work of Spaltenstein \cite{\SPAS},\cite{\SPAF} and on \cite{\HS}. As an 
application we define a natural injective map from the set of unipotent $G$-orbits in $G$ to the set of nilpotent
$G$-orbits in $\fg$ (see 3.3); this maps preserves the dimension of an orbit.

In \cite{\SE}, Serre asked whether a power $u^n$ (where $n$ is an integer not divisible $p$, $p\ge2$) of a 
unipotent element $u\in G$ is conjugate to $u$ under $G$. This is well known to be true when $p\gg0$. In \S2 we 
answer positively this question in general using the theory of Springer's representations; we also discuss an 
analogous property of nilpotent elements.

I wish to thank J.-P. Serre for his interesting questions and comments.

\head 1. Combinatorics\endhead
\subhead 1.1\endsubhead
For $k\in\NN$ let $\ce_k=\{a_*=(a_0,a_1,\do,a_k)\in\NN^{k+1};a_0\le a_1\le\do\le a_k\}$. For $a_*\in\ce_k$ let 
$|a_*|=\sum_ia_i$. For $a_*,a'_*\in\ce_k$ we set $a_*+a'_*=(a_0+a'_0,a_1+a'_1,\do,a_k+a'_k)$. For any $n\in\NN$ 
let $\ce_k^n=\{a_*\in\ce_k;|a_*|=n\}$. We have an imbedding $\ce_k^n@>>>\ce_{k+1}^n$, 
$(a_0,a_1,\do,a_k)\m(0,a_0,a_1,\do,a_k)$. This is a bijection if $k$ is sufficiently large with respect to $n$. 
For $n\in\NN$ let 

$\cc_k^n=\{(a_*,a'_*)\in\ce_k\T\ce_k;|a_*|+|a'_*|=n\}$,

$\cd_k^n=\{(a_*,a'_*)\in\cc_k^n;\text{ either }|a_*|>|a'_*|\text{ or }a_*=a'_*\}$.
\nl
Here $k$ is large (relative to $n$), fixed. Let 

${}^b\cc_k^n=\{(a_*,a'_*)\in\cc_k^n;a'_i\le a_i+2\qua \frl i\in[0,k]\}$,

${}^{b1}\cc_k^n=\{(a_*,a'_*)\in\cc_k^n;a'_i\le a_i+2\qua\frl i\in[0,k],a_i\le a'_{i+1}\qua\frl i\in[0,k-1]\}$,

${}^{b2}\cc_k^n=\{(a_*,a'_*)\in\cc_k^n;a'_i\le a_i+2\qua\frl i\in[0,k],a_i\le a'_{i+1}+2\qua\frl i\in[0,k-1]\}$,

${}^{c1}\cc_k^n=\{(a_*,a'_*)\in\cc_k^n;a_i\le a'_{i+1}+1\qua\frl i\in[0,k-1],a'_i\le a_i+1\qua\frl i\in[0,k]\}$,

${}^d\cd_k^n=\{(a_*,a'_*)\in\cd_k^n;a'_i\le a_i\qua\frl i\in[0,k]\}$,

${}^{d1}\cd_k^n=\{(a_*,a'_*)\in\cd_k^n;a'_i\le a_i\qua\frl i\in[0,k],a_i\le a'_{i+1}+2\qua\frl i\in[0,k-1]\}$,

${}^{d2}\cd_k^n=\{(a_*,a'_*)\in\cd_k^n;a'_i\le a_i\qua\frl i\in[0,k],a_i\le a'_{i+1}+4\qua\frl i\in[0,k-1]\}$.
\nl
Note that 

${}^{b1}\cc_k^n\sub{}^{b2}\cc_k^n\sub{}^b\cc_k^n$,

${}^{c1}\cc_k^n\sub{}^{b2}\cc_k^n\sub\cc_k^n$,

${}^{d1}cd_k^n\sub{}^{d2}cd_k^n\sub{}^d\cd_k^n$.
\nl
The following statements are obvious. If $(a_*,a'_*)\in\cc_k^m$, $(b_*,b'_*)\in\cc_k^{m'}$ then 
$(a_*+b_*,a'_*+b'_*)\in\cc_k^{m+m'}$. If $(a_*,a'_*)\in{}^b\cc_k^{m}$, $(b_*,b'_*)\in{}^d\cd_k^{m'}$, then 
$(a_*+b_*,a'_*+b'_*)\in{}^b\cc_k^{m+m'}$. If $(a_*,a'_*)\in{}^d\cd_k^m$, $(b_*,b'_*)\in{}^d\cd_k^{m'}$ then 
$(a_*+b_*,a'_*+b'_*)\in{}^d\cc_k^{m+m'}$.

In the following result we assume that $k$ is large relative to $n$.

\proclaim{Proposition 1.2} (a) Let $(c_*,c'_*)\in\cc_k^n$. Then either $(c_*,c'_*)\in{}^{c1}\cc_k^n$ or there 
exist $m\ge1,m'\ge1$ such that $m+m'=n$ and $(a_*,a'_*)\in\cc_k^m$, $(b_*,b'_*)\in\cc_k^{m'}$ such that 
$(c_*,c'_*)=(a_*+b_*,a'_*+b'_*)$.

(b) Let $(c_*,c'_*)\in{}^b\cc_k^n$. Then either $(c_*,c'_*)\in{}^{b1}\cc_k^n$ or there exist $m\ge0,m'\ge2$ such 
that $m+m'=n$ and $(a_*,a'_*)\in{}^b\cc_k^m$, $(b_*,b'_*)\in{}^d\cd_k^{m'}$, such that 
$(c_*,c'_*)=(a_*+b_*,a'_*+b'_*)$.

(c) Let $(c_*,c'_*)\in{}^d\cc_k^n$. Then either $(c_*,c'_*)\in{}^{d1}\cc_k^n$ or there exist $m\ge2,m'\ge2$ such 
that $m+m'=n$ and $(a_*,a'_*)\in{}^d\cd_k^m$, $(b_*,b'_*)\in{}^d\cd_k^{m'}$ such that 
$(c_*,c'_*)=(a_*+b_*,a'_*+b'_*)$.
\endproclaim
We prove (a). Assume first that $c_s<c_{s+1}$ for some $s\in[0,k-1]$. Define $(b_*,b'_*)\in\cc^k_r$, $r=k-s>0$, by
$b_i=1$ for $i\in[s+1,k]$, $b_i=0$ for $i\in[0,s]$, $b'_i=0$ for $i\in[0,k]$. Define $(a_*,a'_*)\in\cc^k_{n-r}$ by
$a_i=c_i-1$ for $i\in[s+1,k]$, $a_i=c_i$ in $[0,s]$, $a'_*=c'_*$. We have $a_*+b_*=c_*$, $a'_*+b'_*=c'_*$. If 
$r<n$ we see that (a) holds. If $r=n$ then $(c_*,c'_*)=(b_*,b'_*)\in{}^{c1}\cc_k^n$ and (a) holds again.

Next we assume that $c'_s<c'_{s+1}$ for some $s\in[0,k-1]$. Define $(b_*,b'_*)\in\cc^k_r$, $r=k-s>0$, by $b_i=0$ 
for $i\in[0,k]$, $b'_i=1$ for $i\in[s+1,k]$, $b'_i=0$ for $i\in[0,s]$. Define $(a_*,a'_*)\in\cc^k_{n-r}$ by 
$a_*=c_*$, $a'_i=c'_i-1$ for $i\in[s+1,k]$, $a'_i=c'_i$ for $i\in[0,s]$. We have $a_*+b_*=c_*$, $a'_*+b'_*=c'_*$. 
If $r<n$ we see that (a) holds. If $r=n$ then $(c_*,c'_*)=(b_*,b'_*)\in
{}^{c1}\cc_k^n$ and (a) holds again.

Finally we assume that $c_0=c_1=\do=c_r$, $c'_0=c'_1=\do=c'_r$. Since $k$ is large we can assume that $c_0=0$, 
$c'_0=0$. Then $n=0$ and $(c_*,c'_*)\in{}^{c1}\cc_k^n$.

We prove (b). If $n=0$ we have clearly $(c_*,c'_*)\in{}^{b1}\cc_k^n$. Hence we can assume that $n>0$ and that the 
result is true when $n$ is repaced by $n'\in[0,n-1]$.

Assume first that we can find $0<t\le s\le k$ such that $c'_j=c_j+2$ for $j\in[s+1,k]$, $c'_j<c_j+2$ for 
$j\in[t,s]$, $c_{t-1}<c_t$. Note that if $s<k$ then $c'_s<c'_{s+1}$; indeed, $c'_s<c_s-2\le c_{s+1}-2=c'_{s+1}$.
Define $(b_*,b'_*)\in{}^d\cd^k_r$, $r=2k-t-s+1>0$ by $b_i=1$ for $i\in[t,k]$, $b_i=0$ for $i\in[0,t-1]$, $b'_i=1$
for $i\in[s+1,k]$, $b'_i=0$ for $i\in[0,s]$. Define $(a_*,a'_*)\in{}^b\cc^k_{n-r}$ by $a_i=c_i-1$ for $i\in[t,k]$,
$a_i=c_i$ for $i\in[0,t-1]$, $a'_i=c'_i-1$ for $i\in[s+1,k]$, $a'_i=c'_i$ for $i\in[0,s]$. We have $a_*+b_*=c_*$,
$a'_*+b'_*=c'_*$. If $r\ge2$ we see that (b) holds. If $r=1$ then $t=s=k$ and $a_k=c_k-1$, $a_i=c_i$ for 
$i\in[0,k-1]$, $a'_i=c'_i$ for $i\in[0,k]$. The induction hypothesis is applicable to 
$(a_*,a'_*)\in{}^b\cc^k_{n-1}$. If $(a_*,a'_*)\in{}^{b1}\cc^k_{n-1}$ then clearly 
$(c_*,c'_*)\in{}^{b1}\cc^k_{n-1}$ and (b) holds. If $(a_*,a'_*)\n{}^{b1}\cc^k_{n-1}$ then we can find 
$m\ge0,m'\ge2$ such that $m+m'=n-1$ and $(\ta_*,\ta'_*)\in{}^b\cc_k^m$, $(\tb_*,\tb'_*)\in{}^d\cd_k^{m'}$ such 
that $(a_*,a'_*)=(\ta_*+\tb_*,\ta'_*+\tb'_*)$. Then $(c_*,c'_*)=(\ta_*+\tb_*+b_*,\ta'_*+\tb'_*+b'_*)$ where 
$(\ta_*,\ta'_*)\in{}^b\cc_k^m$, $(\tb_*+b_*,\tb'_*+b'_*)\in{}^d\cd_k^{m'+1}$ so that (b) holds.

Next we assume that $c_i>0$ for some $i$. Then we have $0=c_0=c_1=\do=c_{l-1}<c_l$ for some $l\in[0,k]$. If 
$c'_s<c_s+2$ for some $s\in[l,k]$ then we can assume that $s$ is maximum possible with this property and there are
two possibilities. Either $c'_i<c_i+2$ for all $i\in[l,s]$ and then by the previous paragraph (with $t=l$) we see
that (b) holds; or $c'_i=c_i+2$ for some $i\in[l,s]$ and letting $t-1$ be the largest such $i$ we have $0<t\le s$,
$c'_j<c_j+2$ for $j\in[t,s]$, $c'_j=c_j+2$ for $j=t-1$ and $c_{t-1}=c'_{t-1}-2\le c'_t-2<c_t$; using again the 
previous paragraph we see that (b) holds. Thus we may assume that $c'_i=c_i+2$ for all $i\in[l,k]$. Assume in 
addition that $c'_s<c'_{s+1}$ for some $s\in[l,k-1]$. We can assume that $s$ is maximum possible so that 
$c'_s<c'_{s+1}=\do=c'_k$. We have $c_{s+1}=c'_{s+1}-2>c'_s-2=c_s$ hence $c_s<c_{s+1}$. Define 
$(b_*,b'_*)\in{}^d\cd^k_r$, $r=2k-2s\ge2$, by $b_i=1$ for $i\in[s+1,k]$, $b_i=0$ for $i\in[0,s]$, $b'_i=1$ for 
$i\in[s+1,k]$, $b'_i=0$ for $i\in[0,s]$. Define $(a_*,a'_*)\in{}^b\cc^k_{n-r}$ by $a_i=c_i-1$ for $i\in[s+1,k]$, 
$a_i=c_i$ for $i\in[0,s]$, $a'_i=c'_i-1$ for $i\in[s+1,k]$, $a'_i=c'_i$ for $i\in[0,s]$. We have $a_*+b_*=c_*$,
$a'_*+b'_*=c'_*$. We see that (b) holds. Thus we can assume that $c'_l=c'_{l+1}=\do=c'_k=N+2$ so that 
$c_l=c_{l+1}=\do=c_k=N$. Note that $c'_i\le2$ for $i\in[0,l-1]$. We have $(c_*,c'_*)\in{}^{b1}\cc_k^n$ so that 
(b) holds.

Finally we assume that $c_0=c_1=\do=c_k=0$. Then $c'_i\le2$ for $i\in[0,k]$ and $(c_*,c'_*)\in{}^{b1}\cc_k^n$ so 
that (b) holds. This completes the proof of (b).

We prove (c). If $n=0$ we have clearly $(c_*,c'_*)\in{}^{d1}\cd_k^n$. Hence we can assume that $n>0$ and that the 
result is true when $n$ is repaced by $n'\in[0,n-1]$.

Assume first that we can find $0<t\le s\le k$ such that $c'_j=c_j$ for $j\in[s+1,k]$, $c'_j<c_j$ for $j\in[t,s]$,
$c_{t-1}<c_t$. Note that if $s<k$ then $c'_s<c'_{s+1}$; indeed, $c'_s<c_s\le c_{s+1}=c'_{s+1}$. Define 
$(b_*,b'_*)\in{}^d\cd^k_r$, $r=2k-t-s+1>0$ by $b_i=1$ for $i\in[t,k]$, $b_i=0$ for $i\in[0,t-1]$, $b'_i=1$ for 
$i\in[s+1,k]$, $b'_i=0$ for $i\in[0,s]$. Define $(a_*,a'_*)\in{}^d\cd^k_{n-r}$ by $a_i=c_i-1$ for $i\in[t,k]$, 
$a_i=c_i$ for $i\in[0,t-1]$, $a'_i=c'_i-1$ for $i\in[s+1,k]$, $a'_i=c'_i$ for $i\in[0,s]$. We have 
$a_*+b_*=c_*$, $a'_*+b'_*=c'_*$. If $n-2\ge r\ge2$ we see that (c) holds. If $r=1$ then $t=s=k$ and $a_k=c_k-1$, 
$a_i=c_i$ for $i\in[0,k-1]$, $a'_i=c'_i$ for $i\in[0,k]$. The induction hypothesis is applicable to 
$(a_*,a'_*)\in{}^d\cd^k_{n-1}$. If $(a_*,a'_*)\in{}^{d1}\cd^k_{n-1}$ then clearly 
$(c_*,c'_*)\in{}^{d1}\cd^k_{n-1}$ 
and (c) holds. If $(a_*,a'_*)\n{}^{d1}\cd^k_{n-1}$ then we can find $m\ge2,m'\ge2$ such that $m+m'=n-1$ and 
$(\ta_*,\ta'_*)\in{}^d\cd_k^m$, $(\tb_*,\tb'_*)\in{}^d\cd_k^{m'}$ such that 
$(a_*,a'_*)=(\ta_*+\tb_*,\ta'_*+\tb'_*)$. Then $(c_*,c'_*)=(\ta_*+\tb_*+b_*,\ta'_*+\tb'_*+b'_*)$ where 
$(\ta_*,\ta'_*)\in{}^d\cd_k^m$, $(\tb_*+b_*,\tb'_*+b'_*)\in{}^d\cd_k^{m'+1}$ so that (c) holds. If $r=n-1$ then 
$a_i=0$ for $i\in[0,k-1]$, $a_k=0$, $a'_i=0$ for $i\in[0,k]$; hence $c_i=1$ for $i\in[t,k-1]$, $c_k=2$, $c_i=0$ 
for $i\in[0,t-1]$, $c'_i=1$ for $i\in[s+1,k]$, $c'_i=0$ for $i\in[0,s]$. Hence $(c_*,c'_*)\in{}^d\cd_k^n$ so that
(c) holds. If $r=n$ then $(c_*,c'_*)=(b_*,b'_*)\in{}^d\cd_k^n$ so that (c) holds.

Next we assume that $c_i>0$ for some $i$. Then we have $0=c_0=c_1=\do=c_{l-1}<c_l$ for some $l\in[0,k]$. If 
$c'_s<c_s$ for some $s\in[l,k]$ then we can assume that $s$ is maximum possible with this property and there are 
two possibilities. Either $c'_i<c_i$ for all $i\in[l,s]$ and then by the previous paragraph (with $t=l$) we see 
that (c) holds; or $c'_i=c_i$ for some $i\in[l,s]$ and letting $t-1$ be the largest such $i$ we have $0<t\le s$, 
$c'_j<c_j$ for $j\in[t,s]$, $c'_j=c_j$ for $j=t-1$ and $c_{t-1}=c'_{t-1}\le c'_t<c_t$; using again the previous 
paragraph we see that (c) holds. Thus we may assume that $c'_i=c_i$ for all $i\in[l,k]$. Assume in addition that 
$c'_s<c'_{s+1}$ for some $s\in[l,k-1]$. We can assume that $s$ is maximum possible so that 
$c'_s<c'_{s+1}=\do=c'_k$. We have $c_{s+1}=c'_{s+1}>c'_s=c_s$ hence $c_s<c_{s+1}$. Define 
$(b_*,b'_*)\in{}^d\cd^k_r$ ,$r=2k-2s\ge2$, by $b_i=1$ for $i\in[s+1,k]$, $b_i=0$ for $i\in[0,s]$, $b'_i=1$ for 
$i\in[s+1,k]$, $b'_i=0$ for $i\in[0,s]$. Define $(a_*,a'_*)\in{}^d\cd^k_{n-r}$ by $a_i=c_i-1$ for $i\in[s+1,k]$, 
$a_i=c_i$ for $i\in[0,s]$, $a'_i=c'_i-1$ for $i\in[s+1,k]$, $a'_i=c'_i$ for $i\in[0,s]$. We have $a_*+b_*=c_*$,
$a'_*+b'_*=c'_*$. If $r\le n-2$ we see that (c) holds. If $r=n-1$ then $a_i=0$ for $i\in[0,k-1]$, $a_k=0$, 
$a'_i=0$ for $i\in[0,k]$; hence $c_i=1$ for $i\in[s+1,k-1]$, $c_k=2$, $c_i=0$ for $i\in[0,s]$, $c'_i=1$ for 
$i\in[s+1,k]$, $c'_i=0$ for $i\in[0,s]$. Hence $(c_*,c'_*)\in{}^d\cd_k^n$ so that (c) holds. If $r=n$ then 
$(c_*,c'_*)=(b_*,b'_*)\in{}^d\cd_k^n$ so that (c) holds. Thus we can assume that $c'_l=c'_{l+1}=\do=c'k=N$ so 
that $c_l=c_{l+1}=\do=c_k=N$. Note that $c'_i=0$ for $i\in[0,l-1]$. We have $(c_*,c'_*)\in{}^{d1}\cd_k^n$ so that 
(c) holds.

Finally we assume that $c_0=c_1=\do=c_k=0$. Then $c'_i=0$ for $i\in[0,k]$. In this case we have $n=0$ and
$(c_*,c'_*)\in{}^{d1}\cd_k^n$ so that (c) holds. This completes the proof of (c).

\head 2. On Serre's questions\endhead
\subhead 2.1\endsubhead
For any affine algebraic group $H$ over $\kk$ we denote by $\Lie H$ the Lie algebra of $H$. For any $\co\in\cx_G$
(or $\co\in\cx_\fg$) we set $d_\co=2\dim\cb-\dim\co$. 

\subhead 2.2\endsubhead
We recall the definition of Springer's representations following \cite{\LUG}. Let $\cb$
be the variety of Borel subgroups of $G$. Let $\tcb=\{(g,B)\in G\T\cb;g\in B\}$ and let $f:\tcb@>>>G$ be the first
projection. Let $K=f_!\bbq$. In \cite{\LUG} it was observed that $K$ is an intersection cohomology complex on $G$
coming from a local system on the open dense subset of $G$ consisting on regular semisimple elements. Moreover 
$\WW$ acts naturally on this local system and hence, by "analytic continuation", on $K$. In particular, if 
$\co\in\cx_G$ and $i\in\ZZ$ then $\WW$ acts naturally on the $i$-th cohomology sheaf $\ch^iK|_{\co}$ of 
$K|_{\co}$, an irreducible $G$-equivariant local system on $\co$; hence if $\cl$ is an irreducible $G$-equivariant
local system on $\co$ then $\WW$ acts naturally on the $\bbq$-vector space $\Hom(\cl,\ch^iK|_{\co})$. We denote 
this $\WW$-module (with $i=d_\co$) by $V_{\co,\cl}$. As shown in \cite{\LUI}, $V_{\co,\cl}$ is either $0$ or of 
the form $\bbq\ot E$ where $E\in\Irr(\WW)$; moreover any $E\in\Irr(\WW)$ arises in this way from a unique 
$(\co,\cl)$ and $E\m(\co,\cl)$ is an injective map 

$S_G:\Irr(\WW)@>>>\hcx_G$.
\nl
We would like to define a similar map from $\Irr(\WW)$ to $\hcx_\fg$. Let $\tcb'=\{(x,B)\in\fg\T\cb;x\in\Lie B\}$
and let $f':\tcb'@>>>\fg$ be the first projection. Let $K'=f'_!\bbq$. Now if $p$ is small the set of regular 
semisimple elements in $\fg$ may be empty (this is the case for example if $G=SL_2(\kk)$, $p=2$) so the method of
\cite{\LUI} cannot be used directly. However, T.Xue \cite{\XUE} has observed that the method of \cite{\LUI}, 
\cite{\LUG} can be applied if $G$ is a classical group of adjoint type and $p=2$ (in that case the set of regular
semisimple elements in $\fg$ is open dense in $\fg$). More generally for any $G$ which is adjoint, the set of 
regular semisimple elements in $\fg$ is open dense in $\fg$. (Here is a proof. We 
must only check that if $T$ is a maximal torus of $G$ and $\ft=\Lie T$ then the set $\ft_{reg}$ of regular 
semisimple elements in $\ft$ is open dense in $\ft$. Let $Y=\Hom(\kk^*,T)$. We have $\ft=\kk\ot Y$. Now 
$\ft_{reg}$ is the set of all $x\in\ft$ such that for any root $\a:\ft@>>>\kk$ we have $\a(x)\ne0$. It is enough 
to show that any root $\a:\ft@>>>\kk$ is $\ne0$. We have $\a=1\ot\a_0$ where $\a_0:Y@>>>\ZZ$ is a well defined 
homomorphism. It is enough to show that $\a_0$ is surjective. This follows from the adjointness of $G$.) As in
the group case it now follows that $K'$ is an intersection cohomology complex on $\fg$ coming from a local system
on $\fg_{reg}$. Moreover $\WW$ acts naturally on this local system and hence, by "analytic continuation", on $K'$.
In particular, if $\co\in\cx_\fg$ and $i\in\ZZ$ then $\WW$ acts naturally on the $i$-th cohomology sheaf 
$\ch^iK'|_{\co}$ of $K'|_{\co}$, an irreducible $G$-equivariant local system on $\co$; hence if $\cl$ is an 
irreducible $G$-equivariant local system on $\co$ then $\WW$ acts naturally on the $\bbq$-vector space 
$\Hom(\cl,\ch^iK'|_{\co})$. We denote this $\WW$-module (with $i=d_{\co}$) by $V_{\co,\cl}$. As in \cite{\LUI}, 
\cite{\XUE}, $V_{\co,\cl}$ is either $0$ or of the form $\bbq\ot E$ where $E\in\Irr(\WW)$; moreover any 
$E\in\Irr(\WW)$ arises in this way from a unique $(\co,\cl)$ and $E\m(\co,\cl)$ is an injective map 

$S_\fg:\Irr(\WW)@>>>\hcx_\fg$.
\nl
If $G$ is not assumed to be adjoint, let $G_{ad}$ be the adjoint group of $G$ and let $\fg_{ad}=\Lie G_{ad}$. The
obvious map $\p:\fg@>>>\fg_{ad}$ induces a bijective morphism $\cn_{\fg}@>>>\cn_{\fg_{ad}}$ and a bijection
$\cx_{\fg}@>>>\cx_{\fg_{ad}}$. Now any $G_{ad}$-equivariant irreducible $\bbq$-local system on a $G_{ad}$-orbit in
$\cn_{\fg_{ad}}$ can be viewed as an irreducible $G$-equivariant $\bbq$-local system on the corresponding 
$G$-orbit in $\cn_\fg$. This yields an injective map $\hcx_{\fg_{ad}}@>>>\hcx_\fg$. We define an injective map
$S_\fg:\Irr(\WW)@>>>\hcx_\fg$ as the composition of the last map with $S_{\fg_{ad}}$. 

\subhead 2.3\endsubhead
For any $u\in\cu_G$, let $\cb_u=\{B\in\cb;u\in B\}$ and let $\co$ be the $G$-orbit of $u$ in $\cu_G$. Note that 
$\cb_u$ is a non-empty subvariety of $\cb$ of dimension $d_\co/2$, see \cite{\SPA}. Using this and the definition
of $S_G$ we see that $(\co,\bbq)$ is in the image of $S_G$. Hence there is a well defined injective map 
$S'_G:\cx_G@>>>\Irr(\WW)$ such that for any $\co\in\cx_G$ we have $S'_G(\co)=E$ where $E\in\Irr(\WW)$ is given by
$S_G(E)=(\co,\bbq)$.

Similarly, for any $x\in\cn_\fg$ let $\cb_x=\{B\in\cb;x\in\Lie B\}$ and let $\co$ be the $G$-orbit of $x$ in 
$\cn_\fg$. Note that $\cb_x$ is a non-empty subvariety of $\cb$ of dimension $d_{\co}/2$, see \cite{\HS}. Using 
this and the definition of $S_\fg$ we see that $(\co,\bbq)$ is in the image of $S_\fg$. Hence there is a well 
defined injective map $S'_\fg:\cx_\fg@>>>\Irr(\WW)$ such that for any $\co\in\cx_\fg$ we have $S'_\fg(\co)=E$ 
where $E\in\Irr(\WW)$ is given by $S_\fg(E)=(\co,\bbq)$.

The maps $S'_G,S'_\fg$ can be described directly as follows. For $i\in\ZZ$, we may identify $H^i(\cb)$ ($l$-adic 
cohomology) with the stalk of $\ch^iK$ at $1\in G$ hence the $\WW$-action on $K$ induces a $\WW$-action on the 
vector space $H^i(\cb)$. If $\co\in\cx_G$ and $u\in\co$ then the inclusion $\cb_u@>>>\cb$ induces a linear map 
$f_u:H^{d_\co}(\cb)@>>>H^{d_\co}(\cb_u)$ whose kernel is $\WW$-stable; hence there is an induced action of $\WW$ 
on the image $I_u$ of $f_u$. The $\WW$-module $I_u$ is of the form $\bbq\ot E$ for a well defined $E\in\Irr(\WW)$.
We have $S'_G(\co)=E$. Similarly, if $\co\in\cx_\fg$ and $x\in\co$ then the inclusion $\cb_x@>>>\cb$ induces a 
linear map $\ph_x:H^{d_\co}(\cb)@>>>H^{d_\co}(\cb_x)$ whose kernel is $\WW$-stable; hence there is an induced 
action of $\WW$ on the image $I_x$ of $\ph_x$. The $\WW$-module $I_x$ is of the form $\bbq\ot E$ for a well 
defined $E\in\Irr(\WW)$. We have $S'_\fg(\co)=E$.

Let $\fS_G$ be the image of $S'_G:\cx_G@>>>\Irr(\WW)$. Let $\fS_\fg$ be the image of 
$S'_\fg:\cx_\fg@>>>\Irr(\WW)$.

\subhead 2.4\endsubhead
Any automorphism $a:G@>>>G$ induces a Lie algebra automorphism $a':\fg@>>>\fg$ and an automorphism $\ua$ of $\WW$
as a Coxeter group. Now $a$ (resp. $a'$) induces a permutation $\co\m a(\co)$ (resp. $\co\m a'(\co)$) of $\cx_G$ 
(resp. $\cx_\fg$) denoted again by $a$ (resp. $a'$). Also $\ua$ induces in an obvious way a permutation of 
$\Irr(W)$ denoted again by $\ua$. From the definitions we see that 

$\ua S'_G=S'_Ga$, $\ua S'_\fg=S'_\fg a'$.
\nl
Let $x\m x^p$ be the $p$-th power map $\fg@>>>\fg$ (if $p>1$) and the $0$ map $\fg@>>>\fg$ (if $p=1$). The $r$-th
iteration of this map is denoted by $x\m x^{p^r}$; this restricts to a map $\cn_\fg@>>>\cn_\fg$ which is $0$ for
$r\gg0$. The following result answers questions of Serre \cite{\SE}.

\proclaim{Proposition 2.5}(a) Let $u\in\cu_G$ and let $n\in\ZZ$ be such that $nn'=1$ in $\kk$ for some $n'\in\ZZ$.
Then $u^n$ and $u$ are $G$-conjugate.

(b) Let $x\in\cn_\fg$ and let $x'=a_0x+a_1x^p+a_2x^{p^2}+\do$ where $a_0,a_1,a_2,\do\in\kk$, $a_0\ne0$ (so that 
$x'\in\cn_\fg$). Then $x'$, $x$ are $G$-conjugate.
\endproclaim
We prove (a). Let $\co$ be the $G$-orbit of $u$ and let $\co'$ be the $G$-orbit of $u':=u^n$. Clearly, 
$\cb_u\sub\cb_{u'}$. Since $u'$ is a power of $u$ we have also $\cb_{u'}\sub\cu$ hence $\cb_{u'}=\cb_u$. From 
$\dim\cb_u=\dim\cb_{u'}$ we see that $d_{\co}=d_{\co'}$. The map $f_u:H^{d_\co}(\cb)@>>>H^{d_\co}(\cb_u)$ in 2.3 
remains the same if $u$ is replaced by $u'$. From the description of $S'_G$ given in 2.3 we deduce that 
$S'_G(\co)=S'_G(\co')$. Since $S'_G$ is injective we deduce that $\co=\co'$. This proves (a).

We prove (b). Let $\co$ be the $G$-orbit of $x$ and let $\co'$ be the $G$-orbit of $x'$. Clearly, 
$\cb_x\sub\cb_{x'}$. Since $x=a'_0x'+a'_1x'{}^p+a'_2x'{}^{p^2}+\do$ with $a'_0,a'_1,a'_2,\do\in\kk$, $a'_0=a_0\i$,
we have $\cb_{x'}\sub\cb_x$ hence $\cb_{x'}=\cb_x$. From $\dim\cb_x=\dim\cb_{x'}$ we see that $d_{\co}=d_{\co'}$.
The map $\ph_x:H^{d_\co}(\cb)@>>>H^{d_\co}(\cb_x)$ in 2.3 remains the same if $x$ is replaced by $x'$. From the 
description of $S'_G$ given in 2.3 we deduce that $S'_\fg(\co)=S'_\fg(\co')$. Since $S'_\fg$ is injective we 
deduce that $\co=\co'$. This proves (b).

Parts (a),(b) of the following result answer questions of Serre \cite{\SE}; the proof of (b) below (assuming (a))
is due to Serre \cite{\SE}.

\proclaim{Proposition 2.6}Let $c:G@>>>G$ be an automorphism such that for some maximal torus $T$ of $G$ we have 
$c(t)=t\i$ for all $t\in T$. Let $\tc:\fg@>>>\fg$ be the automorphism of $\fg$ induced by $c$. 

(a) For any $u\in\cu_G$, $c(u),u$ are $G$-conjugate. 

(b) For any $g\in G$, $c(g),g\i$ are $G$-conjugate.

(c) For any $x\in\cn_\fg$, $\tc(x),-x$ are $G$-conjugate. 

(d) For any $x\in\fg$, $\tc(x),-x$ are $G$-conjugate.
\endproclaim
We prove (a). Let $\uc:\WW@>>>\WW$ be the automorphism induced by $c$. If $B\in\cb$ contains $T$ then $T\sub c(B)$
and $B,c(B)$ are in relative position $w_0$, the longest element of $\WW$. Hence if $B,B'$ in $\cb$ contain $T$ 
and are in relative position $w\in\WW$ then $c(B),c(B')$ contain $T$ and are in relative position $w_0ww_0\i$. 
They are also in relative position $\uc(w)$. It follows that $\uc(w)=w_0ww_0\i$ for all $w\in\WW$. Hence the 
induced permutation $\uc:\Irr(\WW)@>>>\Irr(\WW)$ is the identity map. Let $\co$ be the $G$-orbit of $u\in\cu_G$. 
Then $c(\co)$ is the $G$-orbit of $c(u)$. By 2.4 we have $S'_G(c(\co))=\uc(S'_G(\co))=S'_G(\co)$. Since $S'_G$ is
injective it follows that $\co=c(\co)$. This proves (a). 

Following \cite{\SE}, we prove (b) by induction on $\dim(G)$. If $\dim G=0$ the result is trivial. Now assume that
$\dim G>0$. Write $g=su=us$ with $s$ semisimple, $u$ unipotent. If the result holds for $g_1\in G$ then it holds
for any $G$-conjugate of $g_1$. Hence by replacing $g$ by a conjugate we can assume that $s\in T$ so that 
$c(s)=s\i$. Let $Z(s)^0$ be the connected centralizer of $s$, a connected reductive subgroup of $G$ containing 
$T$. Note that $c$ restricts to an automorphism of $Z(s)^0$ of the same type as $c:G@>>>G$. Moreover we have 
$g\in Z(s)^0$. If $Z(s)^0\ne G$ then by the induction hypothesis we see that $c(g),g\i$ are conjugate under 
$Z(s)^0$ hence they are conjugate under $G$. If $Z(s)^0=G$ then by (a), $c(u),u$ are conjugate in $G$. By 2.5(a),
$u,u\i$ are conjugate in $G$. Hence $c(u),u\i$ are conjugate in $G$. In other words, for some $h\in G$ we have 
$c(u)=hu\i h$. Since $s$ is central in $G$ and $c(s)=s\i$ we have $c(s)=hs\i h\i$. It follows that 
$c(g)=c(s)c(u)=hs\i h\i hu\i h=hs\i u\i h\i=hg\i h\i$. This proves (b). 

The proof of (c) is completely similar to that of (a); it uses $S'_\fg$ instead of $S_G$. The proof of (d) is
completely similar to that of (b); it uses (c) and 2.5(b) instead of (b) and 2.5(a).

\head 3. A parametrization of the set of nilpotent $G$-orbits in $\fg$\endhead
\subhead 3.1\endsubhead
Let $V$ be a finite dimensional $\QQ$-vector space. Let $R\sub V^*=\Hom(V,\QQ)$ be a (reduced) root system and let
$W\sub GL(V)$ be the Weyl group of $R$. Let $\Pi$ be a set of simple roots for $R$. Let 
$\Th=\{\b\in R;\b-\a\n R\text{ for all }\a\in\Pi\}$. For any integer $r\ge1$ let $\ca_r$ be the set of all 
$J\sub\Pi\cup\Th$ such that $J$ is linearly independent in $V^*$ and $\sum_{\a\in\Pi}\ZZ\a/\sum_{\b\in J}\ZZ\b$ is
finite of order $r^k$ for some $k\in\NN$. For $J\in\ca_r$ let $W_J$ be the subgroup of $W$ generated by the 
reflections with respect to the roots in $J$. For any $E\in\Irr(W)$ let $b_E$ be the smallest integer $\ge0$ such
that $E$ appears with multiplicity $m_E>0$ in the $b_E$-th symmetric power of $V$ regarded as a $W$-module. Let 
$\Irr(W)^\da=\{E\in\Irr(W);m_E=1\}$. Replacing here $(V,W)$ by $(V,W_J)$ with $J\in\ca_r$ we see that $b_E$ is 
defined for any $E\in\Irr(W_J)$ and that $\Irr(W_J)^\da$ is defined. For $J\in\ca_r$ and $E\in\Irr(W_J)^\da$ there
is a unique $\tE\in\Irr(W)$ such that $\tE$ appears with multiplicity $1$ in $\Ind_{W_J}^WE$ and $b_E=b_{\tE}$; 
moreover, we have $\tE\in\Irr(W)^\da$. We set $\tE=j_{W_J}^WE$. Define $\cs^1_W\sub\Irr(W)^\da$ as in 
\cite{\LUU, 1.3}. Replacing $(V,W)$ by $(V,W_J)$ with $J\in\ca_r$ we obtain a subset 
$\cs^1_{W_J}\sub\Irr(W_J)^\da$. For any integer $r\ge1$ let $\cs^r_W$ be the set of all $E\in\Irr(W)$ such that
$E=j_{W_J}^WE_1$ for some $J\in\ca_r$ and some $E_1\in\cs^1(W_J)$ (see \cite{\LUU, 1.3}). If $r=1$ this agrees 
with the earlier definition of $\cs^1_W$ since in this case $W_J=W$ for any $J\in\ca_r$. For any integer $r\ge1$ 
we define a subset $\ct^r_W$ of $\Irr(W)^\da$ by induction on $|W|$ as follows. If $W=\{1\}$ we set 
$\ct^r_W=\Irr(W)$. If $W\ne\{1\}$ then $\ct^r_W$ is the set of all $E\in\Irr(W)$ such that either $E\in\cs^1_W$ or
$E=j_{W_J}^WE_1$ for some $J\in\ca_r$ with $W_J\ne W$ and some $E_1\in\ct^r(W_J)$. From the definition it is clear
that

$\cs^1_W\sub\cs^r_W\sub\ct^r_W$.
\nl
When $r=1$ we have $\cs^1_W=\ct^r_W$.

We apply these definitions in the case where $r=p$, $V=\QQ\ot\YY_G$ (with $\TT$ being "the maximal torus" of $G$ 
and $\YY_G=\Hom(\kk^*,\TT)$), $R$ is "the root system" of $G$ (a subset of $V^*$) with its canonical set of simple
roots and $W=\WW$ viewed as a subgroup of $GL(V)$. Then the subsets $\cs^1_\WW\sub\cs^p_\WW\sub\ct^p_\WW$ of 
$\Irr(\WW)$ are defined. We can now state the following result.

\proclaim{Proposition 3.2} (a) We have $\fS_G=\cs^p_\WW$.

(b) We have $\fS_\fg=\ct^p_\WW$.
\endproclaim
For (a) see \cite{\LUU, 1.4}. The proof of (b) is given in 3.5. 

\proclaim{Corollary 3.3}There is a unique (injective) map $\t:\cx_G@>>>\cx_\fg$ such that 
$S'_G(\x)=S'_\fg(\t(\x))$ for all $\x\in\cx_G$. 
\endproclaim
The existence and uniqueness of $\t$ follows from $\fS_G\sub\fS_\fg$ which in turn follows from 3.2 and the
inclusion $\cs^p_\WW\sub\ct^p_\WW$. 

It is known that when $p\ne2$ we have $\card\fS_G=\card\fS_\fg$; hence in this case $\t$ is a bijection.

\subhead 3.4\endsubhead
For $n\in\NN$ let $W_n$ be the group of all permutations of the set 

$\{1,2,\do,n,n',\do,2',1'\}$
\nl
which commute 
with the involution $i\m i'$, $i'\m i$; let $W'_n$ be the subgroup of $W_n$ consisting of the even permutations.
Assume that $k\in\NN$ is large relative to $n$. When $G$ is adjoint simple of type $B_n$ or $C_n$ ($n\ge2$) we 
identify $\WW=W_n$ in the standard way; we have a bijection $[a_*,a'_*]\lra(a_*,a'_*)$, 
$\Irr(\WW)=\Irr(W_n)\lra\cc_k^n$ as in \cite{\LCL, 2.3}; moreover, $\Irr(\WW)=\Irr(\WW)^\da$, see 
\cite{\LCL, 2.4}. When $G$ is adjoint simple of type $D_n$ ($n\ge4$) we identify $\WW=W'_n$ in the standard way; 
we have a surjective map $\z:\Irr(\WW)^\da=\Irr(W'_n)^\da@>>>\cd^k_n$ such that for any $\r\in\Irr(W'_n)$ we have
$\z(\r)=(a_*,a'_*)$ where $(a_*,a'_*)\in\cd^k_n$ is such that $\r$ appears in the restriction of $[a_*,a'_*]$ from
$W_n$ to $W'_n$ (the set $\Irr(W'_n)^\da$ is determined by \cite{\LCL, 2.5}); note that $|\z\i(a_*,a'_*)|$ is $2$
if $a_*=a'_*$ and is $1$ otherwise.

\subhead 3.5\endsubhead
In this subsection we prove 3.2(b). We can assume that $G$ is adjoint, simple. If $p=1$ or $p$ is a good prime for
$G$ then $\fS_\fg=\fS_G$ hence using 3.2(a) we have $\fS_\fg=\cs^p_\WW$; in our case we have $\WW_J=\WW$ for any
$J\in\ca_p$ hence from the definitions we have $\cs^p_\WW=\cs^1_\WW=\ct^p_\WW$ and the result follows. In the rest
of this subsection we assume that $p$ is a bad prime for $G$. In this case $\fS_\fg$ has been described explicitly
by Spaltenstein \cite{\SPAS},\cite{\SPAF},\cite{\HS} as follows (assuming that the theory of Springer 
correspondence holds; this assumption can be removed in view of \cite{\XUE} and the remarks in 2.2.) 

If $G$ is of type $C_n$, $n\ge2$ ($p=2$), then we have $\fS_\fg=\Irr(\WW)$. If $G$ is of type $B_n$, $n\ge2$
($p=2$), then, according to \cite{\SPA}, $\fS_\fg=\{[a_*,a'_*]\in\Irr(\WW);(a_*,a'_*)\in{}^b\cc_k^n\}$. (Here $k$
is large and fixed.) If $G$ is of type $D_n$, $n\ge4$ ($p=2$), then $\fS_\fg=\z\i({}^d\cd_k^n)$. If $G$ is of type
$G_2$ ($p=2$ or $3$), of type $F_4$ ($p=3$), of type $E_6$ ($p=2$ or $3$), of type $E_7$ ($p=3$), or of type $E_8$
($p=3$ or $5$) then $\fS_\fg=\fS_G$. If $G$ is of type $F_4$ ($p=2$) then $\fS_\fg=\fS_G\sqc\{1_3,2_3\}$ (notation
as in \cite{\ORA, 4.10}); note that $b_{1_3}=12$, $b_{2_3}=4$). If $G$ is of type $E_7$ ($p=2$) then 
$\fS_\fg=\fS_G\sqc\{84'_a\}$ (notation as in \cite{\ORA, 4.12}; we have $b_{84'_a}=15$). If $G$ is of type $E_8$ 
($p=2$) then $\fS_\fg=\fS_G\sqc\{50_x,700_{xx}\}$ (notation as in \cite{\ORA, 4.13}; we have $b_{50_x}=8$, 
$b_{700_{xx}}=16$). 

On the other hand, for types $B,C,D$, $\ct^2_\WW$ is computed by induction using 1.2, the formulas for the maps
$j_{W_J}^W()$ given in \cite{\LUS, 4.5, 5.3, 6.3} and the known description of $\cs^1_{\WW}$; for exceptional 
types, $\ct^p_\WW$ is computed by induction using the tables in \cite{\ALV} and the known description of 
$\cs^1_{\WW}$. 

In each case, the explicitly described subset $\fS_\fg$ of $\Irr(\WW)$ coincides with the explicitly described
subset $\ct^p_\WW$. This completes the proof of 3.2(b).

To illustrate the inclusion $\fS_\fg\sub\ct^p_{\WW}$ we note that:

if $G$ is of type $E_8$ ($p=2$) then $50_x$, $700_{xx}$ in $\fS_\fg-\fS_G$ are obtained by applying 
$j_{\WW_J}^\WW$ (where $\WW_J$ is of type $E_7\T A_1$) to $15'_a\bxt\sgn$, $84'_a\bxt\sgn$ (which belong to 
$\ct^2_{\WW_J}-\cs^2_{\WW_J}$, $\cs^2_{\WW_J}-\cs^1_{\WW_J}$ respectively);

if $G$ is of type $F_4$ ($p=2$) then $1_3,2_3$ in $\fS_\fg-\fS_G$ are obtained by applying $j_{\WW_J}^\WW$ (where 
$\WW_J$ is of type $B_4,C_3\T A_1$) to an object in $\cs^2_{\WW_J}-\cs^1_{\WW_J}$.

\subhead 3.6\endsubhead
If $G$ is of type $B_n$ or $C_n$, $n\ge2$ ($p=2$), then, according to \cite{\LS}, $\fS_G=\{[a_*,a'_*]\in\Irr(\WW);
(a_*,a'_*)\in{}^{b2}\cc_k^n\}$. (Here $k$ is large and fixed.) If $G$ is of type $D_n$, $n\ge4$ ($p=2$), then 
according to \cite{\LS}, $\fS_G=\z\i({}^{d2}\cd_k^n)$.

\enddocument